\documentclass[11pt,leqno,draft]{article}
\usepackage{amsmath,amssymb,amsthm}

\numberwithin{equation}{section}
\theoremstyle{plain}
\newtheorem{theorem}{Theorem}[section]

\newtheorem{proposition}[theorem]{Proposition}
\theoremstyle{definition}

\theoremstyle{remark}
\newtheorem{remark}{Remark}[section]

\makeatletter
\def\@cite#1#2{[{{\bfseries #1}\if@tempswa , #2\fi}]}
\makeatother

\newcommand{\R}{\mathbb{R}}
\newcommand{\N}{\mathbb{N}}

\begin{document}


\begin{center}
{\Large\bf 
Landau-Lifshitz's conjecture\\
about the motion of a quantum mechanical particle
\\[5pt]under the inverse square potential
}
\end{center}

\begin{center}
~\\[0pt]
Motohiro Sobajima%
\footnote[1]{Department of Mathematics, Tokyo University of Science, 
Dipartimento di Matematica ``Ennio De Giorgi'', Universit\`a del Salento, C.P.193, 73100, Lecce, Italy. E-mail: {\tt msobajima1984@gmail.com}}
and 
Shuji Watanabe%
\footnote[2]{Division of Mathematical Sciences, Graduate School of Engineering, Gunma University,
4-2 Aramaki-machi, Maebashi 371-8510, Japan. E-mail: {\tt shuwatanabe@gunma-u.ac.jp}}

\end{center}

\newenvironment{summary}{\vspace{.5\baselineskip}\begin{list}{}{%
     \setlength{\baselineskip}{0.85\baselineskip}
     \setlength{\topsep}{0pt}
     \setlength{\leftmargin}{15mm}
     \setlength{\rightmargin}{15mm}
     \setlength{\listparindent}{0mm}
     \setlength{\itemindent}{\listparindent}
     \setlength{\parsep}{0pt}
     \item\relax}}{\end{list}\vspace{.5\baselineskip}}
\begin{summary}
{\footnotesize {\bf Abstract.} 
Landau and Lifshitz \cite[Section 35]{LL} conjectured that for an arbitrary $k\in \R$,
there exists the motion of a quantum mechanical particle under the inverse square potential
 $k|x|^{-2}$, $x \in \R^3$.  When $k$ is negative and $| k |$ is very large, 
the inverse square potential becomes very deep and 
generates the very strong attractive force, and hence a quantum mechanical
particle is likely to fall down to the origin (the center of the inverse square potential).
Therefore this conjecture (Landau-Lifshitz's conjecture) seems to be wrong at first sight.
We however prove Landau-Lifshitz's conjecture by showing
that there exists a selfadjoint extension for the Schr\"odinger operator
with the inverse square potential $-\Delta+k|x|^{-2}$ in $\R^N\ (N\geq 2)$
and that the spectrum of the selfadjoint extension is bounded below for an arbitrary $k\in \R$.
We thus give the affirmative and complete answer
to Landau-Lifshitz's conjecture in $\R^N\ (N\geq 2)$. 
} 
\end{summary}



\section{Introduction and main result}

Landau and Lifshitz \cite[Section 35]{LL} conjectured that
for an arbitrary $k\in \R$, there exists the motion of a quantum
mechanical particle under the inverse square potential $k|x|^{-2}$,
$x \in \R^3$. More precisely, for an arbitrary $k\in \R$, Landau and Lifshitz
conjectured that there exists
the motion of such a particle corresponding to the state
of an orbital quantum number $\ell$ satisfying $ \ell (\ell + 1) > -\frac{1}{\, 4\,} -k$.
Here, $\ell\in \N\cup\{0\}$. Consider the case where $k$ is negative and $| k |$ is very large.
The inverse square potential $k|x|^{-2}$ then becomes very deep and 
generates the very strong attractive force, and hence a quantum mechanical
particle is likely to fall down to the origin (the center of the inverse square potential).
In such a case, there does not exist its motion. Therefore, at first sight,
this conjecture (Landau-Lifshitz's conjecture) seems to be wrong.

To give the affirmative and complete answer to Landau-Lifshitz's conjecture
from the viewpoint of operator thoery, one needs to show that for an arbitrary $k\in \R$,
the Schr\"odinger operator $H=-\Delta+k|x|^{-2}$ 
has a selfadjoint extension. Here, $\Delta$ is the Laplacian.  
When a quantum mechanical particle falls down to the origin, 
it is expected that 
the spectrum of the selfadjoint Schr\"odinger operator is not bounded below.
Hence the spectrum of the selfadjoint Schr\"odinger operator should be bounded below
as long as there exists its motion. So one moreover needs to show
that the spectrum of the selfadjoint Schr\"odinger operator
is bounded below. Such a selfadjoint Schr\"odinger operator corresponds to the Hamiltonian
of the physical system. The existence of such a Hamiltonian ensures that
time evolution of the physical system is unitary, and hence ensures
that there exists the motion of a quantum mechanical particle under the inverse square potential.

In this paper we give the affirmative and complete answer
to Landau-Lifshitz's conjecture in $\R^N$, where $N\geq 2$. 
To this end we consider the Schr\"odinger operator 
\begin{equation}\label{obj}
H=-\Delta + \frac{k}{|x|^2}, \quad x\in \R^N
\end{equation}
for an arbitrary $k\in \R$. For such a $k\in \R$, we show that the Schr\"odinger operator 
$H$ defined on a certain set specified later has a selfadjoint extension in $L^2(\R^N)$, \ 
$N\geq 2$ and that the spectrum of the selfadjoint extension is bounded below. 

The Schr\"odinger operator $H$ with $N=1$ appears in the two body problem 
of the Calogero model \cite{Calogero71}, the Calogero-Moser model \cite{moser} 
and the Sutherland model \cite{sutherland}. See also Gitman-Tyutin-Voronov \cite{GTV10} 
and their references. Each model describes a quantum 
mechanical system of many identical particles in one dimension with long-range 
interactions, and has attracted considerable interest. The operator $H$ with $N=1$ 
also appears in Wigner's commutation relations in quantum mechanics (see e.g.
\cite{wigner, yang, OK, OW}), which lead to another quantization in quantum mechanics
called Wigner quantization. In this connection, see also 
\cite{WW, watanabe97, watanabe, watanabe-watanabe, watanabe01}.

First, let $N \geq 3$ and let $k \geq - (N - 2)^2/4$. Let us consider the Schr\"odinger
operator $H$ restricted to 
$C_0^{\infty}(\R^N\setminus\{0\})$ if $N=3, 4$ 
and $C_0^{\infty}(\R^N)$ if $N\geq 5$. 
Then the operator $H$ restricted to 
$C_0^{\infty}(\R^N\setminus\{0\})$ ( or $C_0^{\infty}(\R^N)$)
is nonnegative as a consequence of Hardy's inequality with its optimal constant: 
\begin{equation}\label{o-Hdy}
\frac{(N-2)^2}{4}
\int_{\R^N}\frac{|u|^2}{|x|^2}\,dx
\leq 
\int_{\R^N}|\nabla u|^2\,dx 
\quad {\rm for}\ u\in C_0^\infty(\R^N), \quad N\geq 3.
\end{equation}
By the method of the Friedrichs extension,  the operator $H$
restricted to $C_0^{\infty}(\R^N\setminus\{0\})$ ( or $C_0^{\infty}(\R^N)$)
has a selfadjoint extension and the lower bound
of the spectrum of the selfadjoint extension is zero, and hence the spectrum
of the selfadjoint extension is bounded below if $N \geq 3$ and $k \geq - (N - 2)^2/4$.

Second, let $N \geq 1$ and let $k\geq -(N-2)^2/4+1$.  Let us consider the Schr\"odinger
operator $H$ restricted to $C_0^\infty(\R^N\setminus\{0\})$. Then
Edmunds and Evans \cite[Proposition VII.4.1]{EE} showed
that operator $H$ restricted to $C_0^\infty(\R^N\setminus\{0\})$ 
is essentially selfadjoint if and only if $N \geq 1$ and $k\geq -(N-2)^2/4+1$.
This fact is also stated in Reed and Simon \cite[Theorem X.11]{RS2}.
In this connection, see Okazawa \cite{Okazawa96} for more general potentials.
The lower bound of the spectrum of the closure of  $H$ restricted to
$C_0^\infty(\R^N\setminus\{0\})$ is zero. This is because the result for the case $N = 1, 2, 3, 4$ is
obvious and the result for the case $N \geq 5$ follows from Hardy's inequality.
Thus the operator $H$ restricted to $C_0^\infty(\R^N\setminus\{0\})$ has a 
unique selfadjoint
extension and the spectrum of the selfadjoint extension is bounded below
if $N \geq 1$ and $k\geq -(N-2)^2/4+1$. 

Third, let $N \geq 1$ and let $k >-N/4$. Let us consider the Schr\"odinger operator $H$
restricted to $C_0^{\infty}(\R^N \setminus F)$, where the negligible set
\begin{equation}\label{F}
F = \left\{ x=(x_1,\, x_2,\, \cdots,\, x_N) \in \R^N : x_1 = 0 \;\; \hbox{or} \;\; x_2 = 0 \;\; \hbox{or} \;\; \cdots \;\; \hbox{or} \;\;
x_N =0 \right\}
\end{equation}
is removed from $\R^N$. 
In this case Hardy's inequality \eqref{o-Hdy} cannot be applied if $N\leq 3$.
In spite of this, using a generalized Fourier transform \cite{watanabe02},
Watanabe (the second author of the present paper) showed that
there exists  a Friedrichs extension (selfadjoint extension)
of the operator $H$ restricted to $C_0^\infty(\R^N\setminus F)$
and that the lower bound of the spectrum of the selfadjoint extension is zero.
Therefore, the spectrum of the selfadjoint extension is bounded below 
if $N \geq 1$ and $k>-N/4$.

We are now in a position to state our main result. 

\begin{theorem}\label{thm1}
Suppose $N \geq 2$. For an arbitrary $k\in \R$, let $H$ be in \eqref{obj} and 
let $P\,(\not\equiv 0)$ be an eigenfunction of 
the negative Laplace-Beltrami operator $-\Delta_{S^{N-1}}$ 
corresponding to the eigenvalue $\lambda_\ell=\ell (N-2+\ell)$ for $\ell\in \N\cup\{0\}$ satisfying
\begin{equation}\label{condition}
\lambda_\ell \geq
-\frac{(N-2)^2}{4}-k.
\end{equation}
Define 
\begin{equation}\label{HP}
\begin{cases}
H_Pu:=Hu=-\Delta u+ \dfrac{k}{|x|^2}u, 
\\[8pt]
D(H_P):= C_0^\infty(\R^N\setminus F_P), 
\end{cases}
\end{equation}
where $F_P$ is 
\begin{equation}\label{FP}
F_P:=\{x\in \R^N\setminus\{0\}\;;\; x=\rho\omega,\ \rho>0,\ \omega\in S^{N-1}, \ P(\omega)=0\}\cup \{0\}.
\end{equation}
Then $H_P$ is densely defined in $L^2(\R^N)$ and has a Friedrichs extension.
The lower bound of the spectrum of the Friedrichs extension is zero.
Consequently, for an arbitrary $k\in \R$, there exists the motion of a quantum mechanical particle
under the inverse square potential in $\R^N$, \  $N\geq 2$.
\end{theorem}

\begin{remark}
As mentioned above, the existence of the nonnegative selfadjoint
Schr\"odinger operator $H$ is shown
\begin{itemize}
\item[(1)] 
when $N \geq 3$ and $k \geq - (N - 2)^2/4$,\\[-17pt]
\item[(2)] 
when $N \geq 1$ and $k\geq -(N-2)^2/4+1$, \\[-17pt]
\item[(3)] 
when $N \geq 1$ and $k >-N/4$ or\\[-17pt]
\item[(4)] 
when $N \geq 2$ and $k$ is arbitrary \   (Our theorem (Theorem \ref{thm1})). 
\end{itemize}
Table 1 shows the lower bound of $k$ in each case. Our theorem (Theorem \ref{thm1})
gives the best possible value of $k$ if $N \geq 2$.
\begin{table}[h]\label{fig:values}\caption{The lower bound of $k$}
 \begin{center}
 \begin{tabular}{|l||l|l|l|l|l|l|l|} \hline
\phantom{-0.5} $N$  &\phantom{5} 1 & \phantom{5} 2&
\phantom{5} 3  & \phantom{5} 4  & \phantom{5} 5 & \phantom{5} 6 & $\cdots$ \\ \hline 
$(1) \; -(N - 2)^2/4$ &               &            &-0.25 
 &-1             &-2.25    & -4   & $\cdots$ \\ \hline
$(2) \; -(N-2)^2/4+1$   &\phantom{-}0.75&\phantom{-}1&\phantom{-}0.75
 &\phantom{-}0   &-1.25  & -3    & $\cdots$   \\ \hline
$(3) \; -N/4$         &       -0.25   &-0.5        &-0.75
 &-1             &-1.25    & -1.5   & $\cdots$  \\ \hline
${\rm (4) \  Theorem}$ \ref{thm1} & 
& $-\infty$ & $-\infty$ & $-\infty$
& $-\infty$ & $-\infty$  & $\cdots$ \\ \hline
 \end{tabular}
 \end{center}
\end{table}
\end{remark}

\begin{remark}
If $\ell=N$, that is, $\lambda_N=2N(N-1)$, then we can take 
\begin{equation}\label{Pw}
P(\omega)=x_1 x_2 \cdots x_N,\quad \omega=(x_1, x_2, \ldots, x_N)\in S^{N-1}. 
\end{equation}
In this case, $F_{P}$ coincides with $F$ given by \eqref{F} which is dealt with in \cite{watanabe02}.
Theorem \ref{thm1} asserts that 
$H_P$ with $P$ given by \eqref{Pw} 
has a nonnegative selfadjoint extension 
under the following condition:
\[
k\geq -\left(\frac{3N-2}{2}\right)^2,
\]
which is weaker than that in \cite{watanabe02}.
\end{remark}

This paper is organized as follows. In Section 2, 
we give the Hardy type inequality in 
$D(H_P)$
with optimal constant. 
This is essential in proving the lower spectral bound 
of $H_P$. 
Section 3 is devoted to the proof of Theorem \ref{thm1}.

\section{Hardy type inequality in $D(H_P)$}

In this section we present a Hardy type inequality in $D(H_P)$
with the optimal constant.

\begin{proposition}\label{key}
Let $P$ be an eigenfunction of $-\Delta_{S^{N-1}}$ 
corresponding to the eigenvalue $\lambda_\ell=\ell (N-2+\ell)$ for $\ell\in \N\cup\{0\}$, 
and set $H_P$ be in \eqref{HP}. 
Then for every $u\in D(H_P)$, 
\begin{equation}\label{Hardy}
\left[
\frac{(N-2)^2}{4}
+
\lambda_\ell\right]
\int_{\R^N}\frac{|u|^2}{|x|^2}\,dx
\leq 
\int_{\R^N}|\nabla u|^2\,dx.
\end{equation}
Moreover, the constant in \eqref{Hardy} is optimal.
\end{proposition}

\begin{remark}
We give some comments on the proof of Proposition \ref{key}. 
The auxiliary function $\psi$ in \eqref{psi} plays 
a crucial role in proving the optimality of \eqref{Hardy}.  
In fact, we can observe the optimal constant of \eqref{Hardy} 
in the identity \eqref{c-psi}. 
Furthermore, a family $\{u_m\}$ in \eqref{u-m} approximate to $\psi$ 
works as a minimizing sequence. 
Therefore we conclude that the argument in the proof of Theorem \ref{thm1} 
is natural for deriving \eqref{Hardy} and its optimality.
\end{remark}

\begin{proof}
By the standard approximation argument, 
it suffices to show \eqref{Hardy} for
$u\in C_0^\infty(\R^N\setminus F_{P})$. 
Let $u\in C_0^\infty(\R^N\setminus F_{P})$ and 
put a real-valued function $\psi\in C^\infty(\R^N\setminus\{0\})$ as
\begin{equation}\label{psi}
\psi(x):=\psi(\rho,\omega):= \rho^{-\frac{N}{2}+1}P(\omega). 
\end{equation} 
Then we see from the definition of $F_P$ that 
$\psi(x)\neq 0$ for $x\in \R^N\setminus F_P$ and hence 
$\psi^{-1}\in C^\infty(\R^N\setminus F_P)$. 
Thus setting  
$v=\psi^{-1}u\in C_0^\infty(\R^N\setminus F_{P})$, we have 
\begin{align}
\notag
\int_{\R^N}|\nabla u|^2\,dx
&\,
=
\int_{\R^N}|\nabla( \psi v)|^2\,dx
\\
\notag
&\,
\geq
2\int_{\R^N}\psi\nabla\psi \cdot {\rm Re}(\overline{v}\nabla v)\,dx
+
\int_{\R^N}|\nabla\psi |^2|v|^2\,dx.
\end{align}
Integration by parts gives 
\begin{align}
\notag
2\int_{\R^N}\psi\nabla\psi \cdot {\rm Re}(\overline{v}\nabla v)\,dx
=
-\int_{\R^N}|\nabla \psi|^2|v|^2\,dx
-\int_{\R^N}\psi\Delta \psi|v|^2\,dx.
\end{align}
Combining the above estimates, we see that
\begin{align}
\notag
\int_{\R^N}|\nabla u|^2\,dx
&\,
\geq
-\int_{\R^N}\psi\Delta \psi|v|^2\,dx.
\\
\notag
&\,
\geq
-\int_{\R^N}\psi^{-1}\Delta \psi|u|^2\,dx.
\end{align}
Therefore noting that 
\begin{align}\label{c-psi}
\Delta \psi(x)
&\,=(\rho^{-\frac{N}{2}+1})''P(\omega)
                  +\frac{N-1}{\rho}(\rho^{-\frac{N}{2}+1})'P(\omega)
                  +\frac{1}{\rho^2}(\rho^{-\frac{N}{2}+1}\Delta_{S^{N-1}}P(\omega))
\\
\notag
&\,=
                  -\left(\frac{N-2}{2}\right)^2\rho^{-\frac{N}{2}-1}P(\omega)
                  -\lambda_\ell\rho^{-\frac{N}{2}-1} P(\omega)
\\
\notag
&\,=-\left[\left(\frac{N-2}{2}\right)^2+\lambda_\ell\right]\frac{\psi(x)}{|x|^2}, 
\end{align}
we obtain \eqref{Hardy}:
\begin{align}\label{pr-Hdy}
\int_{\R^N}|\nabla u|^2\,dx
\geq 
\left[\left(\frac{N-2}{2}\right)^2+\lambda_\ell\right]
\int_{\R^N}\frac{|u|^2}{|x|^2}\,dx.
\end{align}

Next we show that $\left(\frac{N-2}{2}\right)^2+\lambda_\ell$ in \eqref{Hardy}
is optimal. 
We may assume without loss of generality that 
$\|P\|_{L^2(S^{N-1})}=1$.
Now we define $C_\ell$ as the optimal constant of \eqref{Hardy}, that is, 
\[
C_P:=\inf\left\{\int_{\R^N}|\nabla u|\,dx
\;;\;u\in D(H_P),\quad 
\int_{\R^N}\frac{|u|^2}{|x|^2}\,dx
=1
\right\}.
\]
Then we see from \eqref{pr-Hdy} that 
\begin{equation}\label{CPLB}
C_P\geq \left(\frac{N-2}{2}\right)^2+\lambda_\ell.
\end{equation}
Conversely, we fix a real-valued function 
$\phi\in C_0^\infty(\R)$ with $\|\phi\|_{L^2(\R)}=1$ and  
choose a family of functions $\{u_m\}_{m\in \N}$ in $C_0^\infty(\R^N\setminus \{0\})$ 
as 
\begin{equation}\label{u-m}
u_{m}(x):=
\frac{1}{m^{\frac{1}{2}}}
\phi\left(\frac{\log |x|}{m}\right)
\psi(x), 
\end{equation}
where $\psi$ is defined in \eqref{psi}.
Then we have $\{u_m\}_{m\in \N}\subset D(H_P)$. 
Using the spherical coordinates and 
change of variables from $\rho$ to $e^{ms}$, 
we obtain 
\begin{align}\label{cul}
\int_{\R^N}\frac{|u_m|^2}{|x|^2}\,dx
&\,=
\frac{1}{m}
   \int_0^\infty
   \frac{1}{\rho}\left|\phi\left(\frac{\log \rho}{m}\right)
   \right|^2
  \,d\rho\times \int_{S^{N-1}}|P(\omega)|^2\,d\omega
\\
\notag
&\,=
\|\phi\|_{L^2(\R)}^2\|P\|_{L^2(S^{N-1})}^2
\\
\notag
&\,=
1,
\end{align}
where we used $\|P\|_{L^2(S^{N-1})}=1$ and $\|\phi\|_{L^2(\R)}=1$. 
On the other hand, 
note that for every $x\in \R^N\setminus\{0\}$, 
\[
\Delta u_m (x)=  
   \frac{1}{m^{\frac{5}{2}}}\phi''\left(\frac{\log |x|}{m}\right)
   -\left[\left(\frac{N-2}{2}\right)^2+\lambda_\ell\right]
   \frac{1}{|x|^2}u_m(x).
\]
Thus using integration by parts (with respect to $x$)
and proceeding the same computation as in \eqref{cul}, 
from the definition of $C_P$ we have 
\begin{align*}
C_P
&\,\leq\int_{\R^N}|\nabla u_m|^2\,dx
\\
&\,=
\int_{\R^N}(-\Delta u_m)\overline{u_m}\,dx
\\
&\,=
-\frac{1}{m^3}
\int_{\R^N}\frac{1}{\rho}
  \phi''\left(\frac{\log\rho}{m}\right)
  \phi  \left(\frac{\log\rho}{m}\right)\,d\rho
\times \int_{S^{N-1}}|P(\omega)|^2\,d\omega
\\
&\,\quad+
\left[\left(\frac{N-2}{2}\right)^2+\lambda_\ell\right]
\int_{\R^N}\frac{|u_m|^2}{|x|^2}\,dx
\\
&\,=
-\frac{1}{m^2}\int_0^\infty\phi''(s)\phi(s)\,ds
+
\left(\frac{N-2}{2}\right)^2+\lambda_\ell.
\end{align*}
Thus integration by parts (with respect to $s$) implies
that for every $m\in \N$,
\begin{equation}\label{CPUB}
C_P\leq 
\frac{1}{m^2}
   \|\phi'\|_{L^2(\R)}^2
   +
   \left(\frac{N-2}{2}\right)^2+\lambda_\ell.
\end{equation}
Therefore 
it follows from \eqref{CPLB} and \eqref{CPUB} that 
$\left(\frac{N-2}{2}\right)^2+\lambda_\ell$ is nothing but 
the best possible constant $C_P$:
\[
C_P=\left(\frac{N-2}{2}\right)^2+\lambda_\ell.
\]
This completes the proof of Proposition \ref{key}.
\end{proof}

\section{Proof of Theorem \ref{thm1}}
\paragraph{}

It is clear that $H_P$ is 
densely defined and  
symmetric in $L^2(\R^N)$. 
Furthermore, 
integration by parts and Proposition \ref{key} imply that 
for every $u\in C_0^\infty(\R^N\setminus F_P)$, 
\begin{align*}
  \int_{\R^N}(H_Pu)\overline{u}\,dx
&\,  =
  \int_{\R^N}|\nabla u|^2\,dx
  +
  k
    \int_{\R^N}\frac{|u|^2}{|x|^2}\,dx
\\
&\,\geq 
  \left[\left(\dfrac{N-2}{2}\right)^2+\lambda_\ell+k\right]
    \int_{\R^N}\frac{|u|^2}{|x|^2}\,dx.
\end{align*}
Therefore from \eqref{condition} we obtain the non-negativity of $H_P$.
Thus 
there exists a Friedrichs extension of $H_P$
(see e.g., Reed-Simon \cite[Theorem X.23]{RS2}). 
This proves Theorem \ref{thm1}. \qed

\bigskip

\noindent \textbf{Acknowledgment}

S. Watanabe is supported in part by the JSPS Grant-in-Aid for Scientific Research (C) 24540112.


{\small

\end{document}